\documentclass[twoside,a4paper,10pt]{article}
\usepackage{amsmath,amsthm,amssymb,latexsym}
\usepackage[latin1]{inputenc}   

\usepackage[dvips]{graphicx}        

\newtheorem{thm}{Theorem}[section]
\newtheorem{theorem}[thm]{Theorem}

\newtheorem{lemma}[thm]{Lemma}

\theoremstyle{remark}

\begin{document}
\title{Hexagonal Tilings and Locally $C_6$ Graphs}

\author{D. Garijo \footnote{ Dep. Matem\'{a}tica Aplicada I. Universidad de
Sevilla. Avda. Reina Mercedes s/n. 41012 Sevilla (Spain).
\{dgarijo,almar,pastora\}@us.es}, I. Gitler\footnote{Dep.
Matemáticas. CINVESTAV-I.P.N. \{igitler@math.cinvestav.mx\}
México D.F.}, A. M\'{a}rquez $^*$, M.P. Revuelta $^*$ }

\date{}

\maketitle

\thispagestyle{empty}

\begin{abstract}
We give a complete classification of hexagonal tilings and locally $C_{6}$ graphs, by showing that
each of them has a natural embedding in the torus or in the Klein bottle (see \cite{thomassen}).
We also show that locally grid graphs, defined in \cite{mmnr, thomassen}, 
are minors of hexagonal
tilings (and by duality of locally $C_6$ graphs) 
by contraction of a particular perfect matching and deletion
of the resulting parallel edges, in a form suitable for 
the study of their Tutte uniqueness.
\end{abstract}

\section{Introduction}

Given a fixed graph $H$, a connected graph $G$ is said to be
locally $H$ if for every vertex $x$ the subgraph induced on the
set of neighbours of $x$ is isomorphic to $H$. For example, if
$P$ is the Petersen graph, then there are three locally $P$
graphs \cite{hall}. In this paper we classify two different
families of graphs, hexagonal tilings and locally $C_6$ graphs.

We first describe all the necessary structures to obtain the classification of hexagonal tilings,
such as the hexagonal cylinder, hexagonal ladder, twisted hexagonal cylinder etc. Some of these
structures appear in \cite{thomassen} in an attempt of classification of these graphs. There exits
an extensive literature on this topic. See for instance the works done by Altshuler
\cite{altshuler1, altshuler2}, Fisk \cite{fisk1, fisk2} and Negami \cite{negami1,negami2}. We also
want to note Ref. \cite{larrion} where locally $C_6$ graphs appear in an unrelated problem. In
this paper, following up the line of research given by Thomassen \cite{thomassen}, we add two new
families to the classification theorem given in \cite{thomassen} and we prove that with these
families we exhaust all the cases. In order to distinguish the families of hexagonal tilings we
study some invariants of graphs such as the chromatic number, shortest essential cycles and
vertex-transitivity. Locally $C_6$ graphs are the dual graphs of hexagonal tilings
\cite{thomassen}, hence the classification theorem of these graphs is obtained from the
classification of hexagonal tilings.

Finally, we are interested in relationships existing between hexagonal tilings, locally $C_6$
graphs and locally grid graphs. Specifically those properties that can be related to different
aspects of the Tutte polynomial. This is a two-variable polynomial $T(G;x,y)$ associated to any
graph $G$, which contains a considerable amount of information about $G$ \cite{bo}. A graph $G$ is
said to be \emph{Tutte unique} if $T(G;x,y)=T(H;x,y)$ implies $G\cong H$ for every other graph
$H$. In \cite{gmr} and \cite{mmnr} the \emph{Tutte uniqueness} of locally grid graphs was studied.
We are interested in a similar study for hexagonal tilings and locally $C_6$ graphs but in a more
unified way, that is in relation to the families of locally graphs that have been Tutte
determined.

Informally a \emph{locally grid graph} is defined as a graph in which the structure around each
vertex is a $3\times 3$ grid (a formal definition is given in Section 3). A complete
classification of these graphs is given in \cite{mmnr,thomassen} and they fall into five families.
Every locally grid graph is a minor of a hexagonal tiling but we are interested in a biyective
minor relationship preserved by duality between hexagonal tilings with the same chromatic number
and locally grid graphs. This specific minor relationship is going to be essential in the study of
the Tutte uniqueness of hexagonal tilings and locally $C_6$ graphs in relation to the Tutte
uniqueness of locally grid graphs. In order to obtain this relation we choose for every family of
hexagonal tilings obtained in the classification theorem, a specific perfect matching, whose
contraction and then deletion of resulting parallel edges (if any) gives rise to a locally grid
graph. There is just one family of hexagonal tilings in which the selected edges are not a
matching. If we select the set of dual edges associated to the perfect matching (hence on the
$C_6$ graph), and we delete them and then contract the set of dual edges associated to the
parallel edges (if any), we obtain the dual of the locally grid graph, which again is a locally
grid graph. These perfect matchings and the set of edges that is not a matching in one of the
families  verify that if we have two hexagonal tilings with the same chromatic number, the results
of the contraction of their perfect matchings (or the set of edges that is not a matching in one
of the families) are two locally grid graphs belonging to different families.

Some standard definitions needed along the paper are: A graph is \emph{$d-$regular} if all
vertices have degree $d$. If $d=3$ the graph is called \emph{cubic}. A \emph{$k-$path} is a graph
with vertices $x_0, x_1, \ldots , x_{k}$ and edges $x_{i-1}x_i$ with $1\leq i\leq k$. A
\emph{$k-$cycle} is obtained from a $(k-1)-$path by adding the edge between the two ends of the
path (vertices of degree one).

\section{Classification of hexagonal tilings and locally $C_6$ graphs}

In this section we give a complete classification of hexagonal tilings, which are connected, cubic
graphs of girth 6, having a collection of $6$-cycles, $\textit{C}$, such that every $2-$path is
contained in precisely one cycle of $\textit{C}$ ($2-$path condition). In particular, a hexagonal
tiling is simple and every vertex belongs to exactly three hexagons (Figure 1). Every hexagon of
the tiling is called a \emph{cell}.

\begin{figure}[htb]
\begin{center}
\includegraphics[width=10mm]{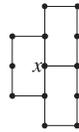} 
\end{center}\caption{Hexagonal structure around $x$}
\label{F0}
\end{figure}

Let $H=P_p\times P_q$ be the $p\times q$ grid, where $P_l$ is a
path with $l$ vertices. Label the vertices of $H$ with the
elements of the abelian group $\mathbb{Z}_p \times \mathbb{Z}_q$
in the natural way. If we add the edges $\{(j,0), (j,q-1)| 0\leq
j \leq p-1 \}$ we obtain a cylinder grid $p\times q$.

A \emph{hexagonal wall of length $k$ and breadth $m$} is defined
as the graph obtained by removing the edges $\{(2i,2j),
(2i+1,2j)\}$ and $\{(2i+1,2j+1), (2i+2,2j+1)\}$ with $0\leq i\leq
\lfloor (m-1)/2 \rfloor$, $0\leq j\leq k-1$ in a $(m+1)\times 2k$
grid. If we delete the same edges in a cylinder grid $(m+1)\times
2k$ the result is a \emph{hexagonal cylinder of length $k$ and
breadth $m$} (Figure 2a). The two cycles of this structure, where
every second vertex has degree two, are called \emph{peripheral
cycles}. Each one of these cycles has $k$ vertices of degree two
labeled as follows: $z_j=(0,2j)$ and $x_j=(m,2j)$ with $0\leq j
\leq k-1$ if $m$ odd, or $z_j=(0,2j)$ and $x_j=(m,2j+1)$ with
$0\leq j \leq k-1$ if $m$ even.

A \emph{hexagonal cylinder circuit of length $k$} is a hexagonal
cylinder of length $k$ and breadth 1. A \emph{hexagonal Möbius
circuit of length $k$} is obtained by adding the edges $\{(0,0),
(1,2k-1)\}$ and $\{(1,0), (0,2k-1)\}$ to a hexagonal wall of
length $k$ and breadth 1. The graph resulting from removing the edges
$\{(0,2j+1), (1,2j+1)| 0\leq j \leq k-1 \}$ and $\{(1,2j),
(2,2j)| 0\leq j \leq k \}$ in a $3\times (2k+1)$ grid, and adding the edges
$\{(0,0),(2,2k)\}$,  $\{(1,0),(1,2k)\}$ and $\{(2,0),(0,2k)\}$ is
called a \emph{parallel hexagonal Möbius circuit} (Figure 2b).

Let $H$ be the $(m+1)\times (2k+m)$ grid. A \emph{hexagonal
ladder of length $k$ and breadth $m$} (Figure 2c) is obtained
by removing the following vertices and edges: $$\begin{tabular}{c}
$\{(j,i)| 0\leq j\leq m-2, 0\leq i\leq m-2-j\}$ \\
\noalign{\smallskip}
$\{(j,2k+i)| 2\leq j\leq m+1, m+1-j\leq i\leq m-1\}$ \\
\noalign{\smallskip} $\{\{(i,m-i+2j),(i+1,m-i+2j)\}| 0\leq i\leq
m-1, 0\leq j\leq  k-1\}$ \end{tabular}$$

\begin{figure}[htb]
\begin{center}
\includegraphics[width=110mm]{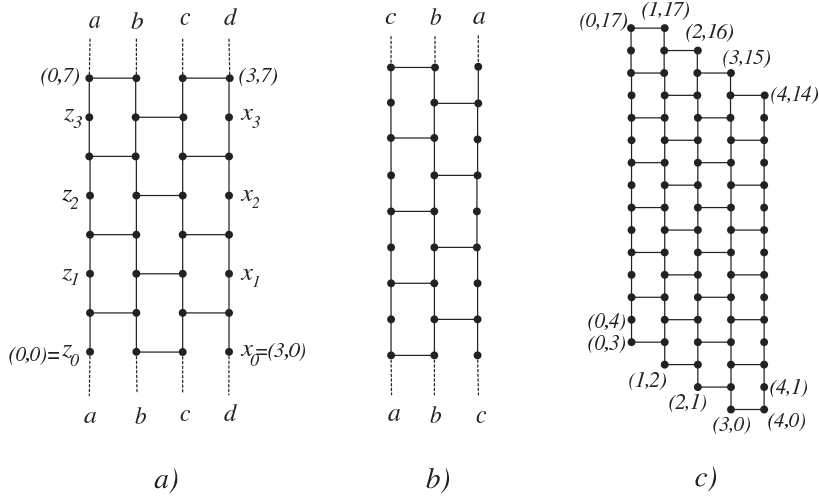} 
\end{center}
\caption{a) Hexagonal cylinder of length 4 and breadth 3  b)
parallel hexagonal Möbius circuit of length 9 c) Hexagonal ladder of
length $7$ and breadth $4$ } \label{CilindromuroescalonMobius}
\end{figure}

From this structure we construct two \emph{twisted hexagonal
cylinders}, $TC_{k,m,1}$ if $k\leq m-2$ (Figure 3b) and
$TC_{k,m,2}$ if $k \geq m+1$ (Figure 3a). The first one is
obtained by adding two vertices, $(0,2k+m)$ and $(m-k-1,3k+2)$, and
the following edges to a hexagonal ladder of length $k$ and
breadth $m$: $\{(0,2k+m),(0,2k+m-1)\}$, $\{(0,2k+m)
,(k+1,m-k-2)\}$, $\{(m-k-1,3k+2),(m-k-1,3k+1)\}$,
$\{(m-k-1,3k+2),(m,0)\}$, $\{(j,2k+m-j),(k+j+1,m-k-j-2)\}$ with
$1\leq j\leq m-k-2$. $TC_{k,m,1}$ also has two peripheral cycles,
$C_1$ and $C_2$, which contain all the vertices of degree two. In
$C_1$, they are $z_j=(0,2k+m-2j)$ and $x_j=(j,m-(j+1))$ with
$0\leq j\leq k$. In $C_2$, $v_0=(m-k-1,3k+2)$,
$v_{i+1}=(m-k+1,3k-i)$, $w_0=(m,2k)$ and $w_{i+1}=(m,2k-(2i+1))$
with $0\leq i \leq k-1$.

To obtain $TC_{k,m,2}$, we delete the vertices $(m+1,i)$, $0\leq
i\leq 2(k-m)-3$ in a hexagonal ladder of length $k$ and breadth
$m+1$, and we add the edges $\{(0,2k+m),(m+1,2(k-m-1)\}$ and
$\{(0,2k+m-(2j+1)),(m,2(k-m-1)-(2j+2))\}$ with $0\leq j\leq
k-m-2$. If $k=m+1$ we do not delete any vertex but we add one
edge, $\{(0,2k+m),(m+1,0)\}$. The vertices of degree two of the
peripheral cycles, $C_1$ and $C_2$, are labeled as follows:
$$\begin{tabular}{rlcll} $C_1$: & $x_i=(m-i,i)$ & and &
$z_i=(0,m+2i+1)$ & with $0\leq i\leq m$.
\\ \noalign{\smallskip}
$C_2$: & $w_i=(i+1,2k+m-i)$ & and & $v_i=(m+1,2k-(2i+1))$ & with
$0\leq i\leq m$. \\ \end{tabular}$$

\begin{figure}[htb]
\begin{center}
\includegraphics[width=100mm]{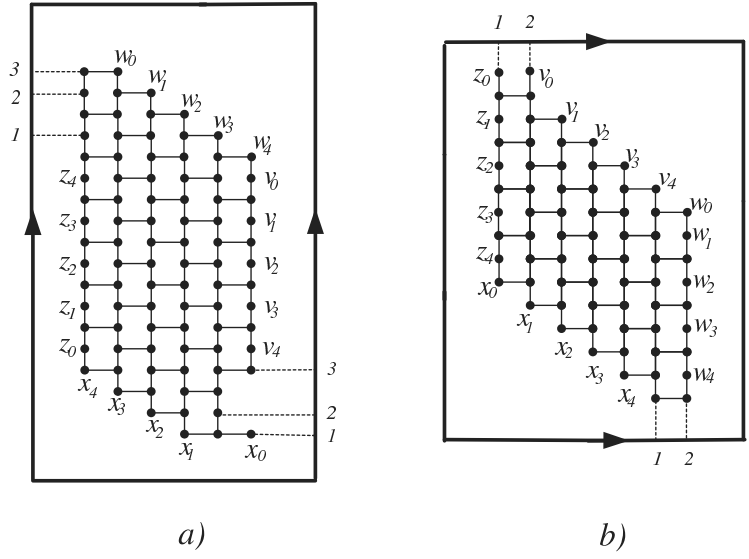} 
\end{center}
\caption{a) $TC_{7,4,2}$ b) $TC_{4,6,1}$  }
\label{cilindrotorcido}
\end{figure}

In order to obtain hexagonal tilings, we must adequately add edges
between the vertices of degree two in the structures defined above.
In the first cases considered below
we add the edges
between vertices on the peripheral cycles
of a hexagonal cylinder. In the last case
we add the edges between vertices on the peripheral cycles of a
twisted hexagonal cylinder.

From a hexagonal cylinder $H$ of length $k$ and
breadth $m$, we obtain the following families of graphs.

\noindent {\bf A)} $H_{k,m,r}$ with $r,k,m\in \mathbb{N}$, \,
$0\leq r\leq \lfloor k/2 \rfloor$,  $m\geq 2$,  $k\geq 3$. If
$m=1$ then $k>3$ and $\lfloor k/2 \rfloor \geq r\geq 2 $ (Figure
5a). $$E(H_{k,m,r})=E(H)\cup \{ \{z_j,x_{j+r} \} | 0\leq j\leq k-1
\}$$

There is a degenerated case, called $H_{k',m',e}$ in
\cite{thomassen}, that we want to stress. It is obtained from a
cycle of even length $k'$ and by adding the adyacencies
$\{z_i,z_{i+m}\}$ taking indices modulo $m$ and taking into
account the $2-$path condition, that is, if $z_i$ is adjacent to
$z_{i+m}$ then $z_{i+1}$ is joined to $z_l$ with $l\equiv (i+1)
({\rm mod} \, m)$ and $l\neq i+m+1$. This graph is a kind of
hexagonal spiral and it is the degenerate case $H_{(k'/2),0,m}$
(Figure 4).

\begin{figure}[htb]
\begin{center}
\includegraphics[width=100mm]{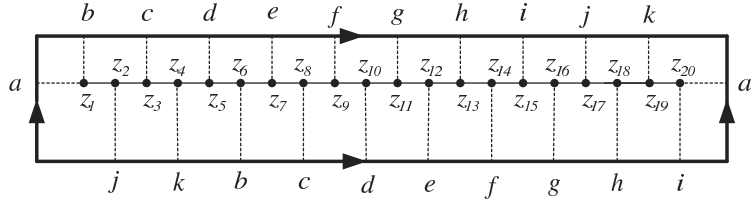} 
\end{center}
\caption{$H_{10,0,5}$} \label{He}
\end{figure}

\noindent {\bf B)} $H_{k,m,a}$ with $m\geq 2$, \, $k\geq 3$
(Figure 5b). $$E(H_{k,m,a})=E(H)\cup \{ \{z_0,x_1\}, \{z_1,x_0\},
\{z_i,x_{k+1-i}\} | 2\leq 1 \leq k-1 \}$$

\noindent {\bf C)} $H_{k,m,b}$ with $k$ even, $m$ odd, \, $m\geq
3$, \, $k\geq 4$ (Figure 5c). $$E(H_{k,m,b})=E(H)\cup \{ \{z_0,x_0
\}, \{z_i,x_{k-1} \} | 1\leq i \leq k-1 \}$$

\begin{figure}[htb]
\begin{center}
\includegraphics[width=130mm]{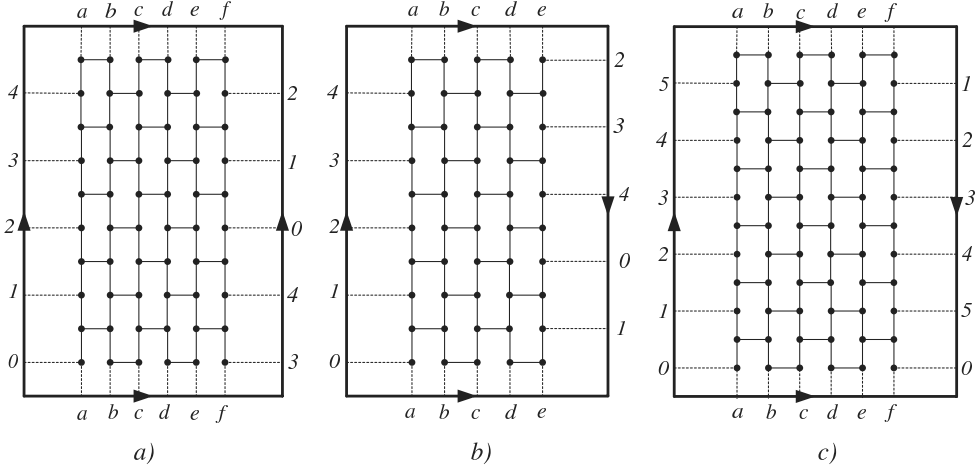} 
\end{center}
\caption{a) $H_{5,5,r}$ b) $H_{5,4,a}$ c) $H_{6,5,b}$}
\label{Hrab}
\end{figure}

\noindent {\bf D)} $H_{k,m,c}$ with $k$ even, \, $k\geq 6$,
$m\geq 1$. (Figure 6). $$E(H_{k,m,c})=E(H)\cup \{ \{z_i,z_{i+k/2}
\}, \{z_i,x_{i+k/2} \} ; \, 0\leq i \leq (k/2)-1 \}$$

An embedding of this graph in the Klein Bottle (Figure 6b) is
obtained by deleting the edges $\{(2i,2j+1),(2i+1,2j+1)\}$ and
$\{(2i+1,2j),(2i+2,2j)\}$ with $0\leq j\leq (k/2)-1$ and $0\leq
i\leq m$ from a $(2m+2)\times k$ grid. Then, we add the edges
$\{(0,2i+1),(2m+1,2i+1)| 0\leq i\leq (k/2)-1\}$ to obtain two
peripheral cycles, whose vertices of degree two are labeled
$z_i=(i,0)$ and $x_i=(i, k-1)$ with $0\leq i\leq 2m+1$. Finally,
we add the edges $\{\{z_i,x_{2m+1-i}\}|0\leq i\leq 2m+1 \}$.

\begin{figure}[htb]
\begin{center}
\includegraphics[width=100mm]{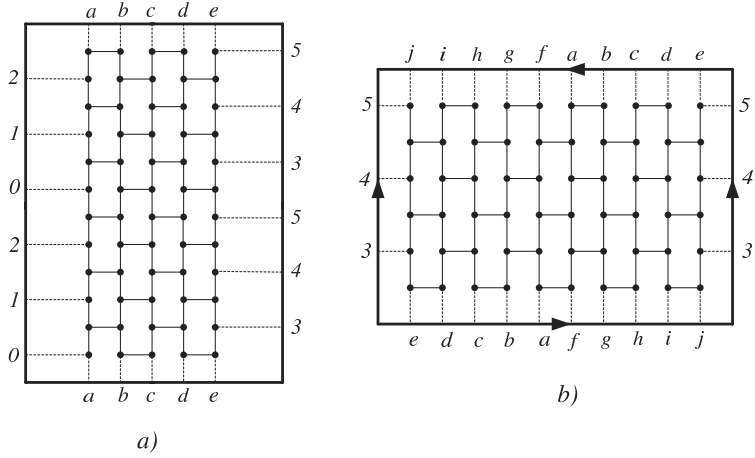} 
\end{center}
\caption{a) $H_{6,4,c}$  b) Embedding of $H_{6,4,c}$ in the Klein
bottle } \label{Hc}
\end{figure}

\noindent {\bf E)} $H_{k,m,f}$ with $k$ odd, $m\geq 0$, $k\geq 7$
(Figure 7).

We add two cycles, $w_0w_1\ldots w_kw_0$ and $v_0v_1\ldots
v_kv_1$, to a hexagonal cylinder of length $k$ and breadth $m$ as
follows: $\{\{z_i,w_{2i}\}, \{ x_i,v_{2i}\} |0\leq i\leq (k-1)/2
\}$. Then, a hexagonal tiling is obtained by adding edges between
the vertices of degree two of this new structure. These edges are
$\{\{z_{(k+2i+1)/2},w_{2i+1}\}, \{ x_{(k+2i+1)/2},v_{2i+1}\}|
0\leq i\leq (k-1)/2 \}$. An embedding of this graph in the Klein
Bottle (Figure 7b), is obtained by deleting the same edges as in
the previous case, from a $(2m+4)\times k$ grid. The edges
$\{(0,2i+1),(2m+3,2i+1)| 0\leq i\leq (k-1/2)-1\}$ are added
giving rise to two peripheral cycles, whose vertices of degree
two are labeled $z_i=(i,0)$ and $x_i=(i, k-1)$ with $0\leq i\leq
2m+3$. Finally, we add the edges $\{\{z_0,x_0\},
\{z_i,x_{2m+4-i}\}|0\leq i\leq 2m+3 \}$.

If $m=0$, we obtain the degenerate case called $H_{k,d}$ in
\cite{thomassen}.

\begin{figure}[htb]
\begin{center}
\includegraphics[width=100mm]{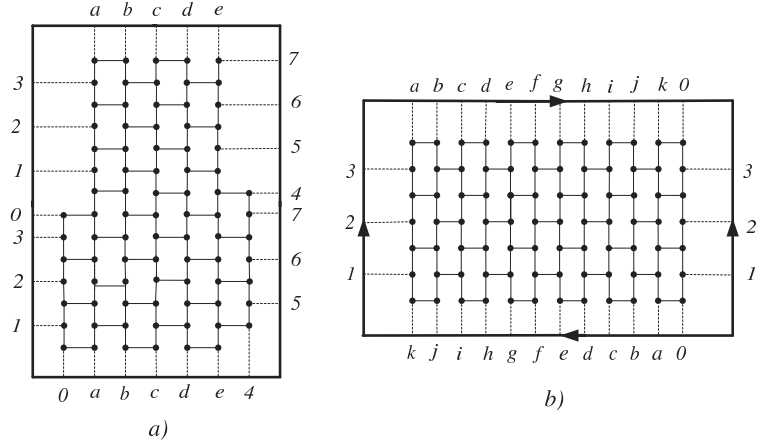} 
\end{center}
\caption{a) $H_{7,4,f}$  b) Embedding of $H_{7,4,f}$ in the Klein
bottle } \label{Hf}
\end{figure}

\noindent {\bf F)} $H_{k,m,g}$ with $k\geq m+1$ and $m\geq 3$
(Figure 8a).

Let $TC_{k,m,2}$ be a twisted hexagonal cylinder of length $k$
and breadth $m$. In order to obtain a hexagonal tiling, we add
the following edges: $$E(H_{k,m,g})=E(TC_{k,m,2})\cup \{ \{z_i,w_i
\}, \{x_i,v_{i} \} ; \, 0\leq i \leq m \}$$

\noindent {\bf G)} $H_{k,m,h}$ with $k\leq m-2$ and $k\geq 2$
(Figure 8b). $$E(H_{k,m,h})=E(TC_{k,m,1})\cup \{ \{z_i,x_i \} |
0\leq i \leq k\} \cup \{v_i,w_{i} \} ; \, 0\leq i \leq k-1 \}$$

\begin{figure}[htb]
\begin{center}
\includegraphics[width=110mm]{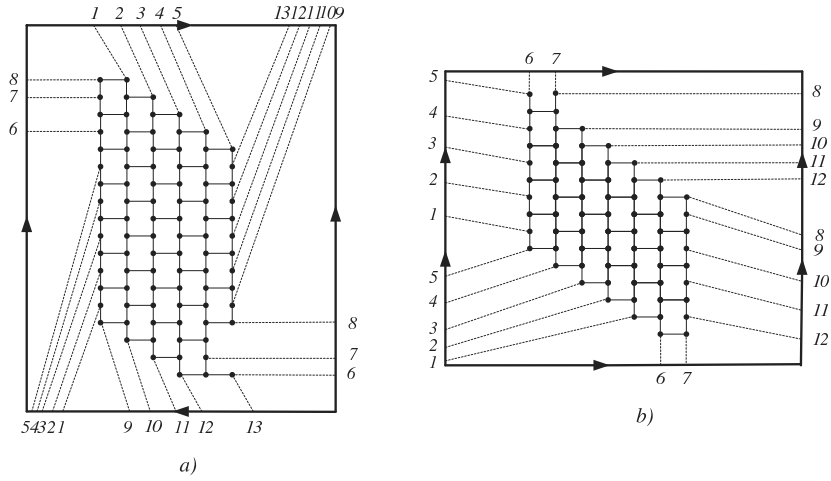} 
\end{center}
\caption{a) $H_{7,4,g}$  b)  $H_{6,4,h}$} \label{Hgh}
\end{figure}

It is straightforward to verify that all the graphs we have
defined are hexagonal tilings. We now prove that these families
exhaust all the possible cases. In order to do so, we study the
shortest essential cycles, the vertex transitivity and the
chromatic number of all the hexagonal tilings defined. Every
hexagonal tiling $G$ has an embedding in the torus or in the Klein
bottle (if $G$ has $v$ vertices, $e$ edges and $h$ hexagons, then
$v=2h$ and $2a=3v$ hence the Euler characteristic is zero).

Given two cycles $C$ and $C'$ in a hexagonal tiling $G$, we say
that $C$ is \emph{locally homotopic} to $C'$ if there exists a
cell, $H$, with $C\cap H $ connected and $C'$ is obtained from $C$
by replacing $C\cap H $ with $H-(C\cap H)$. A \emph{homotopy} is a
sequence of local homotopies. A cycle in $G$ is called
\emph{essential} if it is not homotopic to a cell. This definition
is equivalent to the one given in a graph embedded in a surface
\cite{mmnr}. Let $l_{G}$ be the minimum length of the essential
cycles of $G$, $l_{G}$ is invariant under isomorphism.

\begin{lemma}
Let $G$ be a hexagonal tiling of the types defined in A), B), C),
D), E), F), G) then the length $l_{G}$ of their shortest essential
cycles and the number of these cycles is:

$$\begin{tabular}{|c|c|l|} \hline $G$ & $l_G$ & \\ \hline
$H_{k,m,r}$ & $\begin{tabular}{c} $2k$ \\ \noalign{\medskip}
$2(m+1)$
 \\
\noalign{\medskip} $2(m+1+r-\lfloor (m+1)/2 \rfloor)$ \\
\noalign{\medskip} $2k$ \\ \end{tabular}$ &
$\begin{tabular}{l}  if $k<m+1$\\
\noalign{\medskip}
if $r < \lfloor (m+1)/2 \rfloor < \lfloor k/2 \rfloor$ \\
\noalign{\medskip} if $ \lfloor (m+1)/2 \rfloor \leq r \leq
\lfloor k/2 \rfloor$
\\
\noalign{\medskip} if $k=m+1$ \\
\end{tabular}$ \\ \hline
$H_{k,m,a}$ & min$(2k,2m+2)$ &  \\ \hline $H_{k,m,b}$ &
min$(2k,2m+2)$ & \\ \hline $H_{k,m,c}$ & min$(k+1,4m+4)$ & \\
\hline $H_{k,m,f}$ & min$(k,4m+8)$ & \\  \hline $H_{k,m,h}$ &
$2k+2$ & \\ \hline $H_{k,m,g}$ & $\begin{tabular}{c}
$2(k-m)-2\lfloor (m+1)/2 \rfloor +3$ \\ \noalign{\medskip} $k+2$
 \\
\noalign{\medskip}  $k+3$ \\ \end{tabular}$ &
$\begin{tabular}{l}  if $k>2m+1$\\
\noalign{\medskip}
if $k\leq 2m+1$ and $k$ odd \\
\noalign{\medskip} if $ k<2m+1$ and $k$ even
\end{tabular}$
\\ \hline \end{tabular}$$

$$\begin{tabular}{|c|c|l|} \hline $G$ & number of essential cycles of length $l_G$ &   \\ \hline $H_{k,m,r}$ &  $\begin{tabular}{c} $m+1$  \\
\noalign{\bigskip} $k\left(
\begin{array}{c} m+1 \\ \lfloor (m+1)/2 \rfloor -r \end{array}
\right)$   \\
\noalign{\bigskip} $k\left(
\begin{array}{c} r+\lfloor (m+1)/2 \rfloor \\ m \end{array}
\right)$
\\
\noalign{\bigskip} $m+1+k\left(
\begin{array}{c} m+1 \\ \lfloor (m+1)/2 \rfloor -r \end{array} \right)$ \\
\end{tabular}$ & $\begin{tabular}{l}  if $k<m+1$\\
\noalign{\bigskip} \noalign{\medskip}
if $r <  \lfloor (m+1)/2 \rfloor < \lfloor k/2 \rfloor$ \\
\noalign{\bigskip} \noalign{\medskip}  if $ \lfloor (m+1)/2
\rfloor \leq r \leq \lfloor k/2 \rfloor$
\\
\noalign{\bigskip}  \noalign{\medskip}  if $k=m+1$ \\
\end{tabular}$
 \\ \hline
$H_{k,m,a}$ &  $\begin{tabular}{c} $m+1$ \\
\noalign{\medskip} $2^{m+1}$   \\
\noalign{\medskip} $2^{m+1}+m+1$  \\
\end{tabular}$ & $\begin{tabular}{l} if $k<m+1$ \\
\noalign{\medskip} if $k>m+1$   \\
\noalign{\medskip} if $k=m+1$  \\
\end{tabular}$
 \\  \hline $H_{k,m,b}$ & $\begin{tabular}{c} $m+1$  \\
\noalign{\bigskip} $2\left(
\begin{array}{c} m+1 \\ (m+1)/2 \end{array} \right)+4\sum_{j=1}^{(m-1)/4}\left(
\begin{array}{c} m+1 \\ (m+1)/2 -2j\end{array} \right) $
\\ \noalign{\bigskip}
$m+1+2\left(
\begin{array}{c} m+1 \\ (m+1)/2 \end{array} \right)+4\sum_{j=1}^{(m-1)/4}\left(
\begin{array}{c} m+1 \\ (m+1)/2 -2j\end{array} \right) $  \\
\end{tabular}$ & $\begin{tabular}{l} if $k<m+1$ \\
\noalign{\bigskip} \noalign{\medskip} if $k>m+1$   \\
\noalign{\bigskip} \noalign{\medskip}  if $k=m+1$  \\
\end{tabular}$
\\  \hline $H_{k,m,c}$ & $\begin{tabular}{c}  $(k/2)\left(
\begin{array}{c} 2m+2 \\ m+1 \end{array} \right) $
\\ \noalign{\bigskip}
$2k $  \\
\end{tabular}$ & $\begin{tabular}{l}  if $4m+4<k+1$   \\
\noalign{\medskip} \noalign{\smallskip}  if $4m+4>k+1$  \\
\end{tabular}$
\\  \hline $H_{k,m,f}$ & $\begin{tabular}{c}  $(k-1)/2\left(
\begin{array}{c} 2m+4 \\ m+2 \end{array} \right) $
\\ \noalign{\bigskip}
$2 $  \\
\end{tabular}$ & $\begin{tabular}{l}  if $4m+8<k$   \\
\noalign{\medskip} \noalign{\smallskip}  if $4m+8>k$  \\
\end{tabular}$
\\  \hline   $H_{k,m,h}$ &  $2^{k+1}$ & \\ \hline
$H_{k,m,g}$ & $\begin{tabular}{c}  2
\\ \noalign{\medskip}
$2(k+2)$ \\ \end{tabular}$ & $\begin{tabular}{l}   if $k\leq 2m+1$ and $k$ odd   \\
\noalign{\medskip}   if $k<2m+1$ and $k$ even  \\
\end{tabular}$ \\ \hline
 \end{tabular}$$
\end{lemma}

\

\begin{proof}
We have two different ways of pasting together $j$ ladders each
one containing $i$ hexagons, from which we obtain two structures,
called the \emph{ladder $i\times j$} and the \emph{displaced
ladder $i\times j$ }, shown in Figure 9.

\begin{figure}[htb]
\begin{center}
\includegraphics[width=60mm]{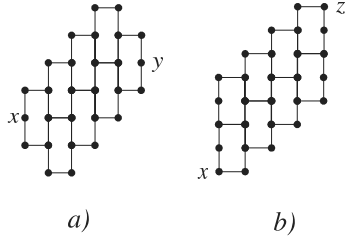} 
\end{center}\caption{ a) Displaced ladder $4\times 2$  b) Ladder $4\times
2$}\label{F12}
\end{figure}

For use below, note that the number of shortest paths between $x$
and $y$ or between $x$ and $z$ in a
ladder $i\times j$ or in a displaced ladder $i\times j$ is $\left(\begin{array}{c} i+j \\
j\end{array}\right)$ and the length of these paths is $2(i+j)-1$.

Recall that every hexagonal tiling defined was obtained by adding
edges to a hexagonal wall or to a hexagonal ladder (except for
$H_{k,m,h}$, in which we also added two vertices). These edges are
called \emph{exterior edges} and every essential cycle must
contain at least one of these edges (for $H_{k,m,h}$ the edges
$\{(0,2k+m), (0,2k+m-1)\}$ and  $\{(m-k-1,3k+2), (m-k-1,3k+1)\}$
are not considered exterior edges).

{\bf (1)} $H_{k,m,r}$

If $k<m+1$, there is only one shortest path determined by each of
the $(m+1)$ exterior edges of the form $\{ (i,0),(i,2k-1)\}$, thus
the resulting cycle has length $2k$.

If $ r < \lfloor (m+1)/2 \rfloor < \lfloor k/2 \rfloor$, the $k$
edges of the form $\{z_i,x_{i+r}\}$ give rise to the shortest
essential cycles. The shortest paths joining the two ends of each
one of these exterior edges have length $2m+1$ and each of them
determine a displaced ladder $r+(m+1)/2 \times (m-1)/2 +1-r$ if
$m$ odd or $(m+2)/2 +r \times (m-2)/2 +1-r$ if $m$ even; hence we
have $k\left(
\begin{array}{c} m+1 \\ \lfloor (m+1)/2 \rfloor -r \end{array}
\right)$ shortest essential cycles of length $2m+2$.

If $\lfloor (m+1)/2 \rfloor \leq r \leq \lfloor k/2 \rfloor$ the
shortest paths in the hexagonal wall that join the two ends of
edges of the form $\{z_i,x_{i+r}\}$ are composed of two parts. The
first part is a path of length $2m+1$ crossing the hexagonal wall
and the second part is a path of length $2r-2\lfloor (m+1)/2
\rfloor$ along a peripheral cycle of the hexagonal cylinder. Each
one of these exterior edges determine a displaced ladder
$r-(m+1)/2 \times m$ if $m$ odd or $r-(m+2)/2  \times m $ if $m$
even.

{\bf (2)} $H_{k,m,a}$

The $m+1$ edges of the form $\{ (i,0),(i,2k-1)\}$ give rise to the
same number of essential cycles as in the previous case.

If $k>m+1$, some exterior edges of the form $\{z_0,x_1 \},
\{z_1,x_0 \}, \{z_i,x_{k+1-i} \}$ with $2\leq i \leq k-1 $
determine shortest essential cycles of length $2m+2$. The shortest
paths that join the two ends of these exterior edges generate
displaced ladders $i\times j$ with $i=m+1-j$ and $0\leq j\leq
m+1$. Hence, the number of shortest essential cycles is
$\sum_{j=0}^{m+1} \left(
\begin{array}{c} m+1 \\ j \end{array}
\right)=2^{m+1}$.

{\bf (3)} $H_{k,m,b}$

If $k<m+1$ we have the same situation as in previous cases. If
$k>m+1$, there are $m+1$ exterior edges of the form $\{z_0,x_0 \},
\{z_i,x_{k-i} \}$ with $1\leq i \leq k-1 $ that determine shortest
essential cycles of length $2m+2$. Two of these edges give rise to
displaced ladders $(m+1)/2\times (m+1)/2$ and the rest of them,
grouped four by four, give rise to displaced ladders
$(m+1)/2+2j\times (m+1)/2-2j$, $0\leq j\leq (m-1)/4$.

\

{\bf (4)} $H_{k,m,c}$

In order to count the shortest essential cycles, we use the
embedding of this graph in the Klein bottle. If $4m+4<k+1$, each
of the $k/2$ exterior edges of the form $\{(0,2i+1),(2m+1,2i+1)
\}$ with $0\leq i \leq (k/2)-1 $ determines a displaced ladder
$(m+1) \times (m+1)$, in which the shortest paths that join the
two ends of this exterior edge is $4m+3$. Hence, the number of
shortest essential cycles is $(k/2)\left(
\begin{array}{c} 2m+2 \\ m+1 \end{array} \right)$ and the length
of these cycles is $4m+4$.

If $4m+4>k+1$, there are just four exterior edges that give rise
to shortest essential cycles of length $k+1$. The shortest paths
that join the two ends of these edges generate ladders
$(k/2)-1\times 1$, therefore the number of shortest essential
cycles is $4\left(
\begin{array}{c} (k/2) \\ 1 \end{array} \right)=2k$.

{\bf (5)} $H_{k,m,f}$

We also use the embedding of this graph in the Klein bottle. If
$4m+8<k$ then reasoning as in the previous case, we obtain
$(k-1)/2$ exterior edges that give rise to displaced ladders
$(m+2) \times (m+2)$. If $4m+8>k$, we just have two shortest
essential cycles of length $k$ generated by the edges $\{
(2m+3,0), (2m+3,k-1)\}$ and $\{ (m-1)/2,0), ((m-1)/2,k-1)\}$.

\,

{\bf (6)} $H_{k,m,h}$

The exterior edge $\{x_0,z_0\}$ determines one shortest essential
cycle of length $2k+2$. The edges $\{(m-k-1,3k+2),(m,0)\}$,
$\{(0,2k+m),(k+1,m-k-2)\}$ and $\{(j,2k+m-j),(k+j+1,m-k-j-2)\}$
with $1\leq j\leq m-k-2$ do not give rise to any shortest
essential cycles. From the remaining, $2k+1$ edges we can
determine different shortest essential cycles, but there are two
exterior edges that generate the same shortest essential cycle.
Every $k+1$ of these edges generate ladders $i\times j$ with
$i=k-j$ and $0\leq j \leq k$ hence the number of shortest
essential cycles is $2\sum_{j=0}^k \left(
\begin{array}{c} k \\ j \end{array}
\right)=2^{k+1}$.

\,

{\bf (6)} $H_{k,m,g}$

If $k=2m+1$ there are just two exterior edges, $\{ z_m,w_m\}$ and
$\{ x_m,v_m\}$, that give rise to shortest essential cycles of
length $2m+3$. If $k<2m+1$ and $k$ even, there are four exterior
edges that determine shortest essential cycles of length $k+3$.
These are the ones that cross the twisted hexagonal cylinder using
$k/2$ hexagons. Each of these exterior edges generate a displaced
ladder $1\times (k/2)$, therefore there are $4(k/2+1)$ shortest
essential cycles. If $k>2m+1$ and $k$ odd, there are two exterior
edges that allow to cross the twisted hexagonal cylinder using
$(k+1)/2$ hexagons and each of these edges give rise to a
displaced ladder $(k+1)/2 \times 0$, hence there are two shortest
essential cycles of length $k+2$. Finally, if $k>2m+1$ the length
of the shortest essential cycles is the sum of the minimum length
of two different paths. The first one, a path that crosses the
hexagonal ladder and the second one, a path in a peripheral cycle
of $TC_{k,m,2}$. In this last case we have not studied the number
of shortest essential cycles.
\end{proof}

From Lemma 2.1, one can prove which of the hexagonal tilings
defined are vertex-transitive graphs and which are not. A graph
$G$ is \emph{vertex-transitive} if for every two vertices of $G$,
$u$ and $v$, there exits an isomorphism of graphs over $V(G)$,
$\sigma $ such that $\sigma (u)=v $. This definition implies that
all the vertices of $G$ have to belong to the same number of
shortest essential cycles.

\begin{lemma}
If $G$ is a hexagonal tiling of the type defined in A), B), C),
D), E), F), G), then $G$ is vertex-transitive if $G$ is isomorphic
to $H_{k,m,r}$ or $H_{4,m,a}$ with $m$ odd.
\end{lemma}

\begin{lemma}
If $G$ is one of the hexagonal tilings defined in A), B), C), D),
E), F), G), then the chromatic number of $G$ is given in the
following table:

$$\begin{tabular}{|c|c|c|c|c|c|c|c|} \hline $G$ & $H_{k,m,r}$ &
$H_{k,m,a}$ & $H_{k,m,b}$ & $H_{k,m,c}$ &
 $H_{k,m,f}$ & $H_{k,m,g}$ & $H_{k,m,h}$ \\ \hline \noalign{\smallskip} \hline
$\chi (G)$ & 2 & 2 & 2 & 3 & 3 & 3 & 2  \\ \hline
\end{tabular} $$
\end{lemma}
\begin{proof}
Let $G$ be one of the hexagonal tilings defined in A), B), C), D),
E), F), G). By Brooks'theorem we know that $\chi(G)<4$ and by
Lemma 2.1 $H_{k,m,c}$, $H_{k,m,f}$ and $H_{k,m,g}$ have cycles of
length odd therefore they can not be bipartite.

It is straightforward to prove that the chromatic number of a
hexagonal cylinder $H$ of length $k$ and breadth $m$ is two for
all $k$ and $m$. Due to the $2-$path condition, the vertices of
degree two of the peripheral cycles have different colors. Hence,
the chromatic number of $H_{k,m,r}$, $H_{k,m,a}$ and $H_{k,m,b}$
is two.

Since every hexagonal ladder admits a $2-$coloring then
 $TC_{k,m,1}$ is bipartite.
We know that $H_{k,m,h}$ is obtained from $TC_{k,m,1}$ by adding
edges between vertices of degree two of the same peripheral cycle.
Now, each peripheral cycle of $TC_{k,m,1}$ has $2k+2$ vertices of
degree two, $k+1$ of these can be assigned the same color and they
are adjacent to the other $k+1$ vertices, which can be assigned
the other color. Therefore the chromatic number of $H_{k,m,h}$ is
two.
\end{proof}

\begin{lemma}
The followings families of hexagonal tilings are not isomorphic:
$$\begin{tabular}{ccll}
 {\bf A)} & $H_{k,m,r}$ & with &  $0\leq r\leq \lfloor k/2 \rfloor $, $m\geq
2$ and $k\geq 3$. If $m=1$ then $k>3$ and $\lfloor k/2 \rfloor
\geq r\geq 2 $. \\ \noalign{\medskip} {\bf B)} & $H_{k,m,a}$ &
with &  $m\geq 2$, $k\geq 3$. \\
\noalign{\medskip} {\bf C)} & $H_{k,m,b}$ & with & $k$ even, $m$
odd, $m\geq 3$, $k\geq 4$. \\ \noalign{\medskip}
 {\bf D)} & $H_{k,m,c}$ & with & $m\geq 1$, $k$ even, $k\geq 6$. \\ \noalign{\medskip}
{\bf E)} & $H_{k,m,f}$ & with & $k$ odd, $m\geq 0$, $k\geq 7$.
\\ \noalign{\medskip} {\bf F)} & $H_{k,m,g}$ &
with & $k\geq m+1$, $m\geq 3$. \\
\noalign{\medskip} {\bf G)} & $H_{k,m,h}$ & with &  $k< m-1$,
$k\geq 2$.
\\ \end{tabular}$$
\end{lemma}
\begin{proof}
By Lemmas 2.2 and 2.3 we just have to prove that the graphs given
in each of the following cases can not be isomorphic.

{\bf (1)} $H_{k,m,a}$ and $H_{k,m,b}$ are not isomorphic since
every graph of the first family contains at most one parallel
hexagonal Möbius circuit and every graph of the second family
contains two.

{\bf (2)} In order to prove that $H_{k,m,a}$ and $H_{k,m,h}$ are
not isomorphic families, we are going to suppose that for every
$k\geq 3$ and $m\geq 2$ there exits $k_1$ and $m_1$ such that
$H_{k,m,a}$ and $H_{k_1,m_1,h}$ are isomorphic and thus obtain a
contradiction. If both graphs are isomorphic, they have the same
number of vertices, shortest essential cycles and the same length
of these cycles, that is, $2k(m+1)=2(k_1+1)(m_1+1)$, $k_1=m$ and
$k=m_1+1$. Hence, the minimum lengths of the non-oriented cycles
of $H_{k,m,a}$ and $H_{k,m,h}$ are $2k$ and $4(m+1)$ respectively,
and thus we reach a contradiction. With an analogous reasoning it
follows that $H_{k,m,b}$ and $H_{k,m,h}$ are not isomorphic
families.

{\bf (3)} In general, the families $H_{k,m,c}$ and $H_{k,m,f}$ can
not have the same number of vertices, and by Lemma 2.1 they cannot
have the same number of shortest essential cycles or the same
length of these cycles.

{\bf (4)} By Lemma 2.1 it is clear that $H_{k,m,g}$ is not
isomorphic neither to $H_{k,m,c}$ nor to $H_{k,m,f}$ because the
length and the number of shortest essential cycles do not coincide
in these graphs.
\end{proof}

\newpage

\begin{theorem}
If $G$ is a hexagonal tiling with $N$ vertices, then one and only
one of the following holds:

$$\begin{tabular}{lcl}
 {\bf A)} & $G\simeq H_{k,m,r}$ & with  $N=2k(m+1)$, $0\leq r\leq \lfloor k/2 \rfloor $, $m\geq
2$, $k\geq 3$. If $m=1$ then $k>3$  \\ \noalign{\smallskip} & &
and $\lfloor k/2 \rfloor \geq r\geq 2 $.   \\ \noalign{\medskip}
{\bf B)} & $G\simeq H_{k,m,a}$ &
with  $N=2k(m+1)$, $m\geq 2$, $k\geq 3$. \\
\noalign{\medskip} {\bf C)} & $G\simeq H_{k,m,b}$ & with
$N=2k(m+1)$, $k$ even, $m$ odd, $m\geq 3$, $k\geq 4$. \\
\noalign{\medskip}
 {\bf D)} & $G\simeq H_{k,m,c}$ & with  $N=2k(m+1)$, $m\geq 1$, $k$ even, $k\geq 6$. \\ \noalign{\medskip}
{\bf E)} & $G\simeq H_{k,m,f}$ & with  $N=2k(m+2)$, $k$ odd,
$m\geq 0$, $k\geq 7$. \\ \noalign{\medskip} {\bf F)} & $G\simeq
H_{k,m,g}$ &
with  $N=2(m+1)(k+2)$, $k\geq m+1$, $m\geq 3$. \\
\noalign{\medskip} {\bf G)} & $G\simeq H_{k,m,h}$ & with
$N=2(m+1)(k+1)$, $k< m-1$, $k\geq 2$.
\\ \end{tabular}$$
\end{theorem}

\begin{proof}
The argument of the proof is essentially the same as the one given in \cite{thomassen}. The
difference between both proofs is that we include two new families to the list given in Theorem
3.1 of \cite{thomassen}, $H_{k,m,g}$ and $H_{k,m,h}$. We consider the families $H_{k,d}$ and
$H_{k,m,e}$ of \cite{thomassen} as degenerated cases of the families $H_{k,m,f}$ and $H_{k,m,r}$
respectively. Therefore we just study the case in which $G$ is a hexagonal tiling containing a
hexagonal cylinder circuit of length $k$. We can extend this circuit either to a hexagonal
cylinder of length $k$ and maximum breadth $m$, or to one of the two twisted hexagonal cylinders,
$TC_{(k/2)-1,m,1}$ or $TC_{l,(k/2)-1,2}$. The first case is studied in \cite{thomassen} obtaining
the families $H_{k,m,a}$, $H_{k,m,b}$, $H_{k,m,r}$, $H_{k,m,c}$ and $H_{k,m,f}$.

Assume that the hexagonal cylinder circuit is extended to a
twisted hexagonal cylinder whose peripheral cycles, $C_1$ and
$C_2$ are labeled as shown in Figure 3. If some vertex of
$C_1$ is joined to some vertex of $C_2$, then by the $2-$path
condition every vertex of degree two of $C_1$ has to be joined to
a vertex of degree two of $C_2$. In a twisted hexagonal cylinder,
we have that each two vertices of degree two of a peripheral
cycle are at distance two except the couples $x_0z_k$, $w_0w_1$
and $v_0v_1$ in $TC_{k,m,1}$ and $z_0x_k$, $w_kv_0$ in
$TC_{k,m,2}$. These couples determine the forms of joining vertices
of degree two in order to obtain a hexagonal tiling. There are two
possibilities, if $z_i$ is joined to $v_i$ and $x_i$ to $w_i$, we
are in the case studied in \cite{thomassen} in which we
extend the circuit to a hexagonal cylinder. If $z_i$ is joined to
$w_i$ and $x_i$ to $v_i$, $G$ is isomorphic to $H_{l,(k/2)-1,g}$.

Assume now that no vertex of $C_1$ is adjacent to a vertex of
$C_2$. Every vertex of degree two of each peripheral cycle has to
be joined to another vertex of degree two of the same peripheral
cycle. There is just one possibility, $z_i$ joined to $x_i$ and
$v_i$ to $w_i$, thus $G$ is isomorphic to $H_{k/2-1,m,h}$.
\end{proof}

The \emph{geometric dual} graph $G^*$ of a graph
$G$ is a graph whose vertex set is formed by the faces of $G$ and
two vertices are adjacent if the corresponding faces share an edge.

\begin{theorem}
\cite{thomassen} Let $G'$ be a connected $6$-regular graph and
$C'$ a collection of $3-$cycles in $G'$ such that, for every
vertex $v$ of $G'$, there are precisely six cycles in $C'$ that
contain $v$ and their union is a $6-$wheel $W_{v}$ with $v$ as
center. Suppose further that $G'$ has no nonplanar subgraphs of
radius 1. Then $G'$ is a dual graph of a hexagonal tiling.
\end{theorem}

Theorems 2.5 and 2.6 give us a complete classification of locally
$C_6$ graphs.

\section{Relation with Locally Grid Graphs}

In this section we establish a biyective minor relationship preserved by duality between
hexagonal tilings with the same chromatic number and locally grid graphs. We want to remark that
this minor relationship between hexagonal tilings, locally $C_6$ graphs and locally grid graphs is
essential in the study of the Tutte uniqueness. The locally grid condition is different in that it
involves not only a vertex and its neighbours, but also four vertices at distance two. Let $N(x)$
be the set of neighbours of a vertex $x$. We say that a $4-$regular, connected graph $G$ is a
\emph{locally grid graph} if for every vertex $x$ there exists an ordering $x_1,x_2,x_3,x_4$ of
$N(x)$ and four different vertices $y_1,y_2,y_3,y_4$, such that, taking the indices modulo 4,
$$\begin{tabular}{ccl}
$N(x_i)\cap N(x_{i+1}) $ & $=$ & $\{x,y_i\}$ \\
\noalign{\medskip} $N(x_i)\cap  N(x_{i+2}) $ & $=$ & $\{x\}$ \\
\end{tabular}$$ \noindent and there are no more adjacencies among $\{
x, x_1, \ldots ,x_4,y_1, \ldots , y_4\}$ than those required by
this condition (Figure 10).

\begin{figure}[htb]
\begin{center}
\includegraphics[width=10mm]{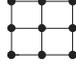} 
\end{center}
\caption{Locally Grid Structure}\label{grid0}
\end{figure}

A locally grid graph is simple, two-connected, triangle-free, and
every vertex belongs to exactly four squares (cycles of length 4).
A complete classification of
locally grid graphs appears in \cite{mmnr}. They fall into
several families and each of them has a natural embedding in the
torus or in the Klein bottle.

Let $H=P_p\times P_q$ be the $p\times q$ grid, where $P_l$ is a
path with $l$ vertices. Label the vertices of $H$ with the
elements of the abelian group $\mathbb{Z}_p \times \mathbb{Z}_q$
in the natural way.

{\bf The Torus} $T_{p,q}^{\delta}$  with $ \, p\geq 5$, $ \, 0\leq
\delta\leq p/2$, $ \, \delta + q\geq 5 \, $ if $ \, q\geq 4$, $ \,
\delta + q\geq 6 \, $ if $ \, q=2,3 \, $ or $ \, 4\leq \delta <
p/2 \, $ with $ \, \delta \neq p/3,p/4 \, $ if $ \, q=1$.
$$\begin{tabular}{cclcl} $E(T_{p,q}^{\delta})$ & $=$ & $E(H)$ & $
\cup$ & $ \{
\{(i,0),(i+\delta,q-1)\}, \, 0\leq i\leq p-1 \}$ \\
\noalign{\medskip} & & & $\cup$ & $ \{ \{(0,j),(p-1,j)\}, \, 0\leq j\leq q-1 \}$ \\
\end{tabular}$$

{\bf The Klein Bottle} $K_{p,q}^1$ with $ \, p\geq 5$, $ p$ odd, $
 \, q\geq 5$.
$$\begin{tabular}{cclcl} $E(K_{p,q}^1)$ & $=$ & $E(H)$ & $ \cup$
& $ \{
\{(j,0),(p-j-1,q-1)\}, \, 0\leq j\leq p-1 \}$ \\
\noalign{\medskip} & & & $\cup$ & $ \{
\{(0,j),(p-1,j)\}, \, 0\leq j\leq q-1 \}$ \\
\end{tabular}$$

{\bf The Klein Bottle} $K_{p,q}^0$ with $p\geq 5$, $p$ even,
 $q\geq 4$.
$$\begin{tabular}{cclcl}  $E(K_{p,q}^0)$ & $=$ & $E(H)$ & $ \cup$
& $ \{
\{(j,0),(p-j-1,q-1)\}, \, 0\leq j\leq p-1 \}$ \\
\noalign{\medskip} & & & $\cup$ & $ \{
\{(0,j),(p-1,j)\}, \, 0\leq j\leq q-1 \}$ \\
\end{tabular}$$

{\bf The Klein Bottle} $K_{p,q}^2$ with $p\geq 5$, $p$ even,
 $q\geq 5$.
$$\begin{tabular}{cclcl} $E(K_{p,q}^2)$ & $=$ & $E(H)$ & $ \cup$
& $ \{
\{(j,0),(p-j,q-1)\}, \, 0\leq j\leq p-1 \}$ \\
\noalign{\medskip} & & & $\cup$ & $ \{
\{(0,j),(p-1,j)\}, \, 0\leq j\leq q-1 \}$ \\
\end{tabular}$$

{\bf The graphs $S_{p,q}$} with $p\geq 3$ and $q\geq 6$.

\,

$ \begin{tabular}{lcclcl} If $ p\leq q$ & $E(S_{p,q})$ & $=$ &
$E(H)$ & $ \cup $ &
$\{(j,0),(p-j,q-p+j)\}, \, 0\leq j\leq p-1 \}$ \\
\noalign{\medskip} & & & & $\cup$ & $ \{ \{(0,i),(i,q-1)\}, \,
0\leq i\leq p-1 \}$ \\ \noalign{\medskip} & & & & $ \cup$ & $ \{
\{(0,i),(p-1,i-p)\}, \, p\leq i\leq q-1 \}$\\
\end{tabular}$

\

\,

\,

$\begin{tabular}{lcclcl} If $q\leq p$ &  $E(S_{p,q})$ & $=$ &
$E(H)$ & $ \cup$ & $ \{
\{(j,0),(0,q-1-j)\}, \, 0\leq j\leq q-1 \}$ \\
\noalign{\medskip} & & & & $\cup$ & $ \{ \{(p-1-i,q-1),(p-1,i)\},
\, 0\leq i\leq q-1 \}$ \\ \noalign{\medskip} & & & & $\cup$ & $ \{
\{(i,q-1),(i+q,0)\}, \, 0\leq i\leq p-q-1 \}$\\
\end{tabular}$

\newpage

\begin{lemma}
If $G$ is a locally grid graph then $G^*=G$ if $G\in \{T_{p,q}^r,
K_{p,q}^1, S_{p,q}\}$ and $(K_{p,q}^0)^*=K_{p,q}^2$.
\end{lemma}

\begin{proof}
Let $H$ be the $p\times q$ grid. $H$ has $(p-1)(q-1)$ squares (cycles of length four). If we
replace every square by a vertex and two vertices are adjacent if the corresponding squares share
an edge, then the resulting graph is a $(p-1)\times (q-1)$ grid. To construct locally grid graphs,
we add edges between vertices of degree two and three of $H$. That is, we add $p+q-1$ squares.
Now, if $G$ is a locally grid graph with $pq$ vertices, then $G^*$ has $pq$ vertices and it is
obtained by adding edges between vertices of degree two and three of a $p\times q$ grid, denoted
$H'$. Vertices of $H$ and $H'$ are labeled $(i,j)$ and $(i,j)^*$ respectively, for $0\leq i\leq
p-1$ and $0\leq j\leq q-1$. Due to the classification theorem of locally grid graphs \cite{mmnr},
we can consider the following cases.

{\bf (1)} If $G\simeq T_{p,q}^{\delta}$, every vertex $(0,j)^*$
is associated to the square with vertices $(0,j-1)$, $(1,j-1)$,
$(0,j+1)$ and $(1,j+1)$ then it has to be adjacent to the vertex
of the square $(0,j-1)$, $(p-1,j-1)$, $(0,j+1)$, $(p-1,j+1)$,
that is, $(p-1,j)^*$. Now, vertices $(i,0)^*$ and
$(i+\delta,q-1)^*$ have to be adjacent since the squares $(i,0)$,
$(i+1,0)$, $(i+\delta,q-1)$, $(i+\delta+1,q-1)$ and
$(i+\delta,q-1)$, $(i+\delta+1,q-1)$, $(i+\delta,q-2)$,
$(i+\delta+1,q-2)$ share an edge. Hence $V(G^*)\simeq V(G)$ and
$E(G^*) \simeq E(G)$. (Figure 11a)

The cases in which $G\simeq K_{p,q}^{1}$ or $G\simeq S_{p,q}$ are
similar to case (1) and we omit the proof for sake
of brevity.

{\bf (2)} If $G\simeq K_{p,q}^{0}$, reasoning as in (1) every
vertex $(0,j)^*$ is adjacent to $(p-1,j)^*$. The vertices
$(p-1,0)^*$ and $(p-1,q-1)^*$ have to be adjacent since the
squares $(p-1,0)$, $(0,0)$, $(p-1,q-1)$, $(0,q-1)$ and
$(p-1,q-1)$, $(p-1,q-2)$, $(0,q-2)$, $(0,q-1)$ share an edge.
Since $p$ is even, $G^*$ is a locally grid graph and since it has
two adjacencies, by \cite{mmnr} $G^*$ it is isomorphic to
$K_{p,q}^{2}$. (Figure 11b)
\end{proof}

\begin{figure}[htb]
\begin{center}
\includegraphics[width=150mm]{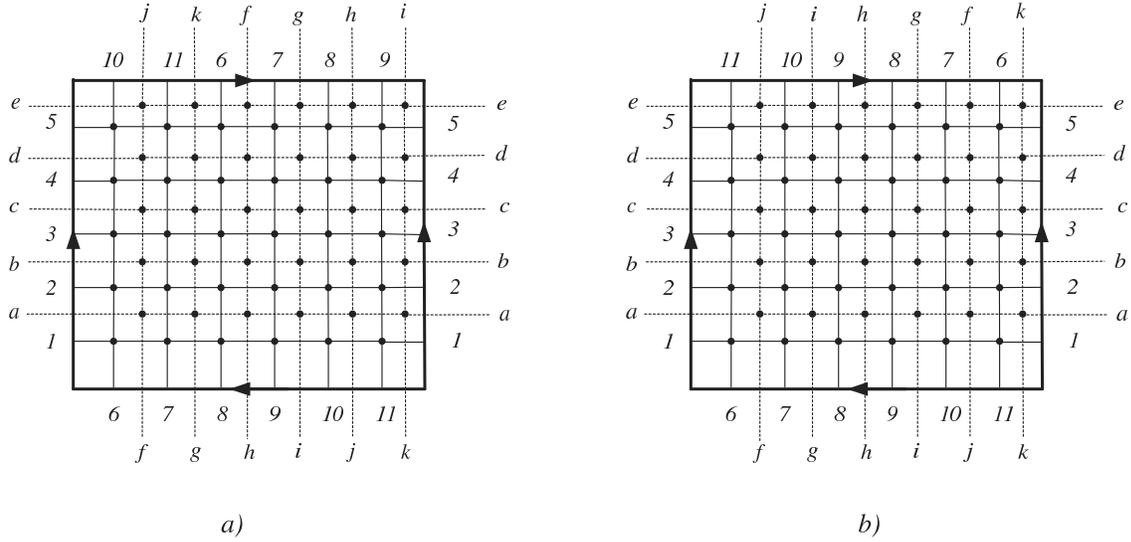} 
\end{center}
\caption{ a) $T^2_{6,5}$ and its dual, $T^2_{6,5}$ (dotted) b)
$K^0_{6,5}$ and its dual, $K^2_{6,5}$ (dotted)
 }\label{dualesgrid}
\end{figure}

\newpage

\begin{theorem}
Locally grid graphs are minors of hexagonal tilings (and by duality of locally $C_6$ graphs) by
contraction of a perfect matching and deletion of the resulting parallel edges, except in one case
in which a set of edges that do not form matching is contracted.

\end{theorem}

\begin{proof}
Let $H$ be a hexagonal tiling, by Theorem 2.5 we know that $H$ belongs to one of the families
$H_{k,m,r}$, $H_{k,m,a}$, $H_{k,m,b}$, $H_{k,m,c}$, $H_{k,m,f}$, $H_{k,m,g}$ and $H_{k,m,h}$. In
order to prove that locally grid graphs are minors of hexagonal tilings, we are going to select a
perfect matching in each one of the families, except in $H_{k,m,f}$ in which the selected edge set
is not a matching. Then we obtain the locally grid graph by contracting the edges of this
matching, and deleting parallel edges if necessary. There are just two cases in which we have to
delete parallel edges, $H_{k,m,f}$ and $H_{k,m,h}$.

Let $C$ be a hexagonal cylinder of length $k$ and breadth $m$. We
can take the following perfect matching, $P$, in $C$: $\{\{
(i,2j), (i,2j+1)\} | 0\leq i \leq m, 0\leq j\leq k-1\}$. Then by
contracting the edges in $P$ we obtain a $k\times (m+1)$ cylinder
grid. This matching is also a perfect matching of $H_{k,m,r}$,
$H_{k,m,a}$ and $H_{k,m,b}$ and no exterior edge of these graphs
is contained in $P$. Therefore, we obtain by contracting $P$ in
$H_{k,m,r}$ the graph $T_{k,m+1}^r$, in $H_{k,m,a}$ we obtain the
graph $K_{k,m+1}^0$ if $k$ even and $K_{k,m+1}^1$ if $k$ odd, and
in $H_{k,m,b}$ we obtain the graph $K_{k,m+1}^2$.

To select a perfect matching in $H_{k,m,c}$, we use the embedding
of this graph in the Klein bottle and we take the same
perfect matching that was specified in the previous case. By
contracting the edges of $P$ we obtain the graph $K_{2m+2, k/2}^0$.

The case of $H_{k,m,f}$ is slightly different. We take the
embedding of this graph in the Klein bottle, and consider the
hexagons with vertices $(0,l)(0,l+1)(0,l+2)(1,l)(1,l+1)(1,l+2)$
where $l=2, \ldots , 2k-5$. Then, $P$ is given by $P_1\cup P_2$,
with: $$\begin{tabular}{c}
 $P_1=\{\{ (i,2j), (i,2j+1)\} | 0\leq i \leq m, 0\leq j\leq
l/2\}$ \\ \noalign{\medskip} $P_2=\{ \{(i,l+1),(i,l+2)\}, \ldots
\{(i,2k-2),(i,2k-1)\}| 0\leq i \leq m\}$ \end{tabular}$$ If the
edges of $P$ are contracted and we delete the resulting parallel edges, we
obtain $K_{2m+4,(k-1)/2}^2$.

For an illustration of these operations see the example given in
Figure 13. In this example we start from $H_{7,4,f}$ and $P$ is
given by the dotted edges. After contracting the selected edge set
and deleting the resulting parallel edges we obtain $K^2_{12,3}$.

For a hexagonal
ladder of length $k$ and breadth $m$, we take the following edges for the matching:
 $\{(0,m-1),(0,m)\}$,
$\{(0,m+1),(0,m+2)\}$,$ \ldots $ , $\{(0,2k+m-3),(0,2k+m-2)\}$,
\linebreak $\{(i,m-(i+1)),(i,m-i)\}$,
$\{(i,m-(i+1)+2),(i,m-i+2)\}$, $ \ldots $ ,
$\{(i,2k+m-(i+1))$,\linebreak $(i,2k+m-i)\}$ with $1\leq i\leq m$.
In order to obtain a perfect matching in $H_{k,m,h}$ we add the
edges $\{(0,2k+m-1), (0,2k+m)\}$ and  $\{(m-k-1,3k+2), (m,0)\}$.
Then by contracting the edges of this matching and deleting the
resulting parallel edge we obtain the locally grid graph
$S_{m+1,k+1}$ with $k<m$. If we consider a similar selection of
edges in a hexagonal ladder of length $k$ and breadth $m+1$ and
we add $\{(0,2k+m),(m+1,2(k-m-1))\}$ we obtain a perfect matching
of $H_{k,m,g}$ whose contraction gives the graph $S_{m+1,k+2}$
with $k\geq m+1$ (Figure 14).

For locally $C_6$ graphs, we follow the same procedure as for
hexagonal tilings. In each case we take the set of dual edges
associated to $P$ (Figure 12).

\begin{figure}[htb]
\begin{center}
\includegraphics[width=30mm]{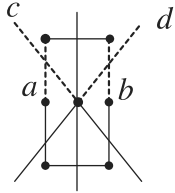} 
\end{center}
\caption{$c$ and $d$ are the dual edges of $a$ and $b$
respectively}\label{aristadual}
\end{figure}

If $G$ is the locally grid graph obtained from the contraction of
the edges of $P$ and deletion of the resulting parallel edges in a
hexagonal tiling $H$, then $G^*$ is obtained by the deletion of
the set $P^*$ of dual edges associated to the perfect matching $P$
and contraction of the set of dual edges associated to the
resulting parallel edges in a locally $C_6$ graph $H^*$. By
Theorems 2.5 and 2.6 and Lemma 3.1, all the cases are determined.
Figures 13 and 14 show two examples. In Figure 13, we start from
$H_{7,4,f}^*$ selecting the dual edges of those belonging to the
selected edge set of $H_{7,4,f}$. After applying the minor
operations we obtain $K^0_{12,3}$, that is the dual graph of
$K^2_{12,3}$. In Figure 14, we delete the dual edges of those
belonging to the perfect matching of $H_{7,4,g}$ obtaining
$S_{4,9}$.

To conclude:
$$\begin{tabular}{rcl}
\begin{tabular}{|c|c|} \hline \begin{tabular}{c} Hexagonal \\
tiling \\ \end{tabular} & \begin{tabular}{c} Minor by contraction \\ and deletion of \\
parallel edges \\ \end{tabular} \\
\hline $H_{k,m,r}$ &  $T_{k,m+1}^r$  \\ \hline
 $H_{k,m,a}$ &
$\begin{tabular}{ccc} $K_{k,m+1}^0$ & if & $k$ even \\
\noalign{\smallskip} $K_{k,m+1}^1$ & if & $k$ odd \\
\end{tabular}$
\\ \hline $H_{k,m,b}$ & $K^2_{k,m+1}$  \\
 \hline $H_{k,m,c}$ & $K_{2m+2,k/2}^0$  \\
\hline $H_{k,m,f}$ & $K_{2m+4,(k-1)/2}^2$
\\ \hline $H_{k,m,g}$ & $S_{m+1,k+2}$ \\
 \hline $H_{k,m,h}$ & $S_{m+1,k+1}$ \\ \hline \end{tabular}
& $\overleftrightarrow{dual}$  & \begin{tabular}{|c|c|} \hline \begin{tabular}{c} Locally \\
$C_6$ graph \\ \end{tabular} & \begin{tabular}{c} Minor by deletion and \\  contraction of dual edges \\
of parallel edges \\ \end{tabular} \\
\hline $H_{k,m,r}^*$ &  $T_{k,m+1}^r$ \\ \hline $H_{k,m,a}^*$ &
$\begin{tabular}{ccc} $K_{k,m+1}^2$ & if & $k$ even \\
\noalign{\smallskip} $K_{k,m+1}^1$ & if & $k$ odd \\
\end{tabular}$
\\ \hline $H_{k,m,b}^*$ & $K^0_{k,m+1}$  \\
\hline $H_{k,m,c}^*$ & $K_{2m+2,k/2}^2$  \\
\hline $H_{k,m,f}^*$ & $K_{2m+4,(k-1)/2}^0$
\\ \hline $H_{k,m,g}^*$ & $S_{m+1,k+2}$ \\
 \hline $H_{k,m,h}^*$ & $S_{m+1,k+1}$ \\ \hline \end{tabular}
\end{tabular}$$
\end{proof}

\begin{figure}[htb]
\begin{center}
\includegraphics[width=130mm]{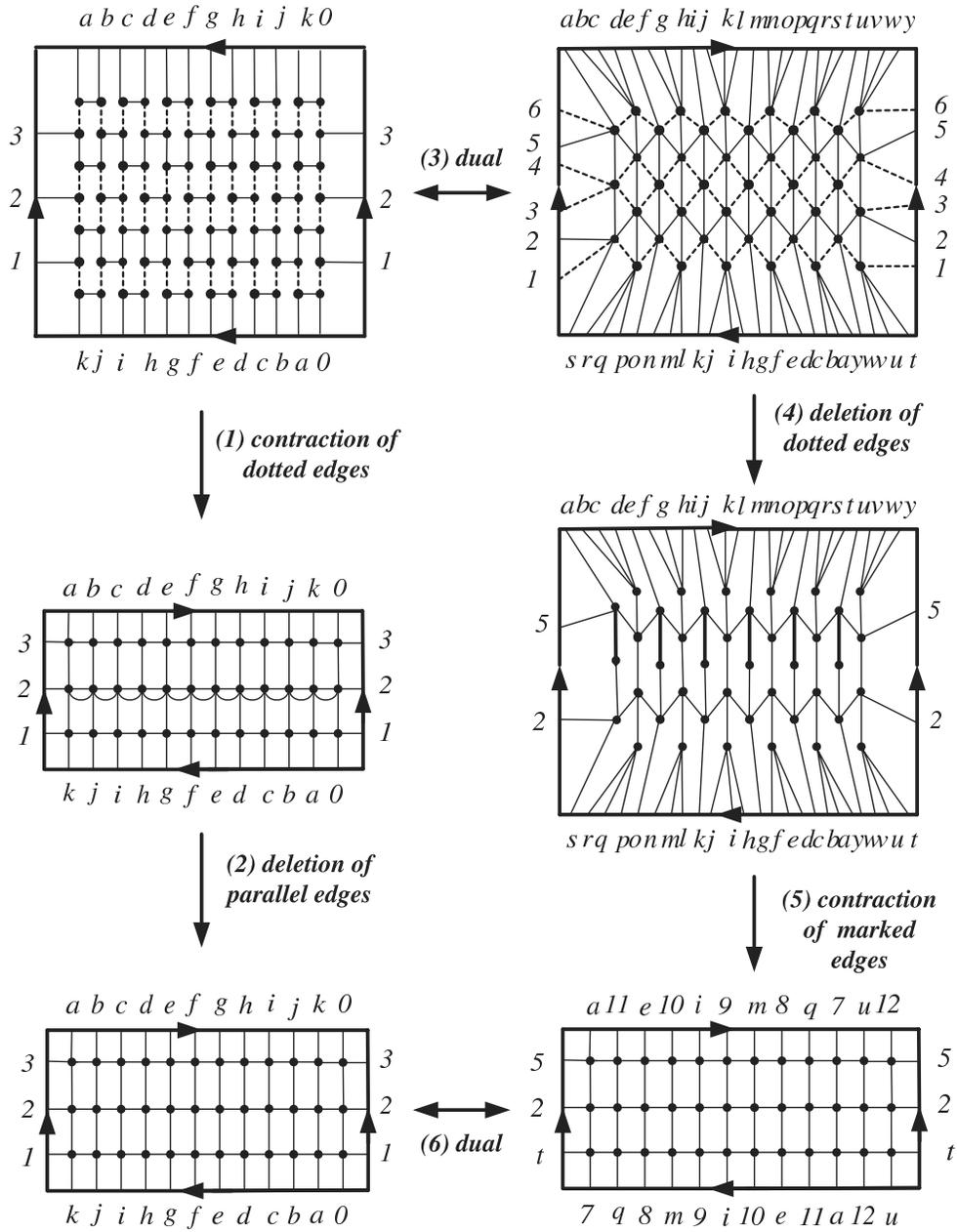} 
\end{center}
\caption{ Deletion and contraction of the selected edge set in
$H_{7,4,f}$ and $H_{7,4,f}^*$ }\label{perfectmatching}
\end{figure}

\begin{figure}[htb]
\begin{center}
\includegraphics[width=150mm]{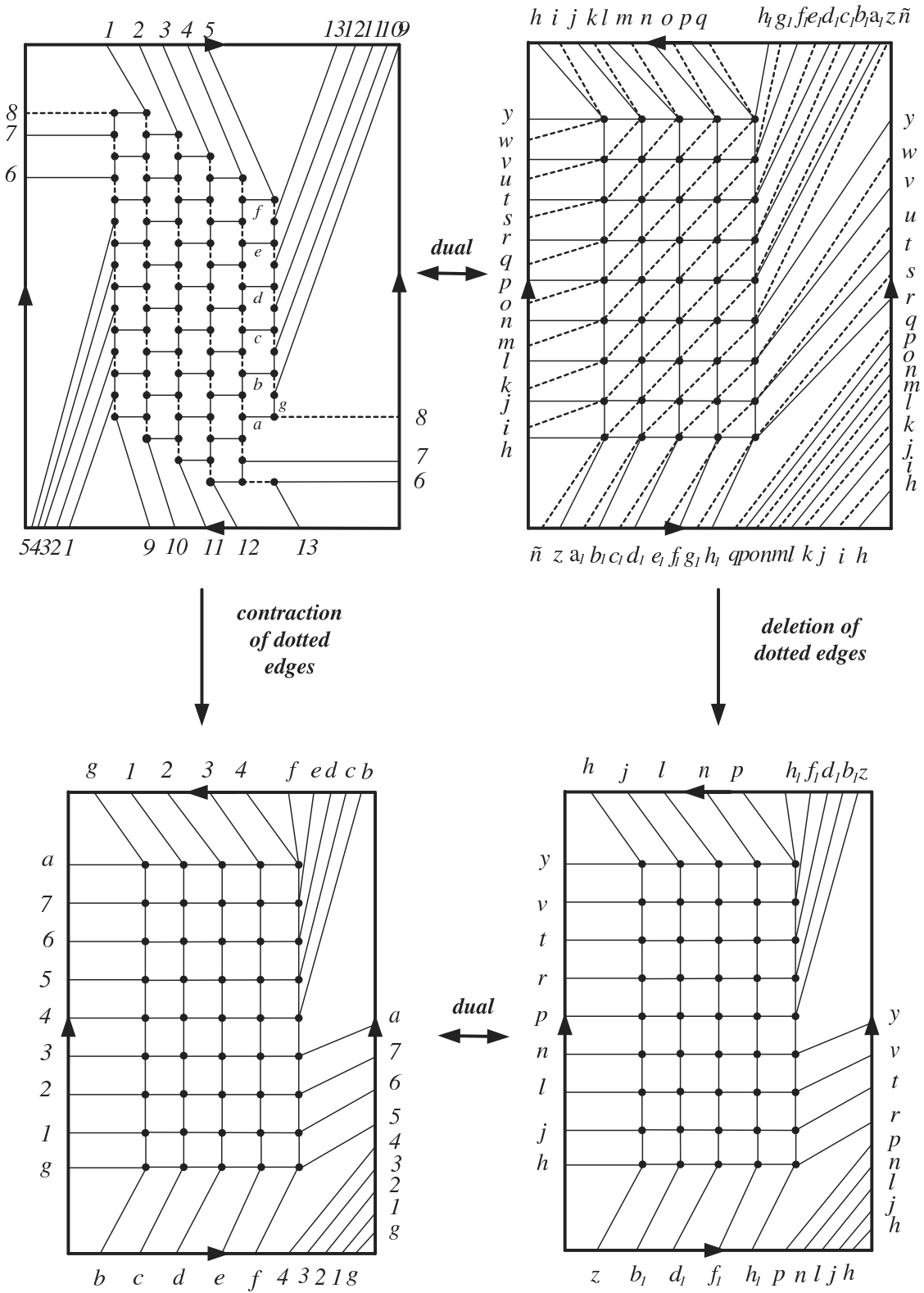} 
\end{center}
\caption{ Deletion and contraction of the edges of a perfect
matching in $H_{7,4,g}$ and
$H_{7,4,g}^*$}\label{matchingHg}\end{figure}

\end{document}